\documentclass[12pt]{article}
\usepackage[log-declarations=false]{xparse}

\usepackage[colorlinks=true,citecolor=black,linkcolor=black,urlcolor=blue]{hyperref}
\usepackage{float}

\usepackage{tikz}

\usepackage[abbrev,logical-quotes,msc-links]{amsrefs}

\usepackage{amssymb}
\usepackage{amsthm}

\newtheorem{theorem}{Theorem}

\newtheorem{lemma}[theorem]{Lemma}
\newtheorem{corollary}[theorem]{Corollary}

\newtheorem{conjecture}[]{Conjecture}

\theoremstyle{plain}

\newcommand{\stack}[6]{$$ \mbox{ #1 } \mbox{ #2 } \left\{
	   \begin{array}{ll}
	   \mbox{ #3 } & \mbox{ #4 } \\ \\
	   \mbox{ #5 } & \mbox{ #6 }
	   \end{array} \right. $$ }

\newcommand{\nmsystem}{$(n,m)$-structure}		

\usepackage{enumitem}
\setlist[enumerate,1]{label=(\roman*)}

 {\begin{list}{}%
         {\setlength{\leftmargin}{#1}}%
         \item[]%
 }
 {\end{list}}

\title{On a conjecture of Szemer\'edi and Petruska}

\author{Adam S. Jobson \qquad Andr\'{e} E. K\'{e}zdy \qquad Tim Pervenecki \\
	\small Department of Mathematics\\[-0.8ex]
	\small University of Louisville\\[-0.8ex]
	\small Louisville, Kentucky, U.S.A.\\
	\small\tt kezdy@louisville.edu}

\begin{document}
\maketitle

\begin{abstract}
Consider a $3$-uniform hypergraph of order $n$ with
clique number $k$ such that the intersection of all its $k$-cliques is empty. Szemer\'edi and Petruska \cite{SzemerediPetruska}  proved  $n\leq 8m^2+3m$, for fixed $m=n-k$, and they conjectured the sharp bound $n\leq{m+2\choose 2}$. Tuza \cite{Tuza2} proved the best known bound, $n\leq \frac{3}{4}m^2+m+1$, using the machinery of $\tau$-critical hypergraphs.
Here we propose an alternative approach, combining a decomposition process introduced by Szemer\'edi and Petruska with the skew version of Bollob\'as's theorem to prove $n\leq m^2 + 6m + 2$.  While the bound obtained here is weaker than Tuza's bound,
it is a proof-of-concept for a different approach and a call to apply dimension bounds from linear algebra.
\end{abstract}

\section{Introduction}
\label{Introduction}

Let ${\cal N} = \{N_1,\ldots,N_\ell\}$ be a collection of $k$-subsets of $[n] = \{1,\ldots,n\}$.
Set $V = \bigcup_{i=1}^\ell N_i$.
Assume that $n = |V|$, $\ell \geq 2$, and $k \geq 3$.  Set $m = n - k$; that is, $\left|\overline{N_i}\right| =|V - N_i|= m$.
We further assume that ${\cal N}$ satisfies the following two properties:
\begin{quote}
\begin{itemize}
\item[(i)] $\bigcap_{i=1}^\ell N_i = \varnothing$, but $\bigcap_{j\neq i} N_j \neq \varnothing$ for all $i=1,\ldots,\ell.$
\item[(ii)] For any $X \subseteq V$ such that $|X| = k+1$, there exists a subset $T \subseteq X$ such that $|T| = 3$ and 
$T \not\subseteq N_i$, for all $i=1,\ldots,\ell$.
\end{itemize}
\end{quote}

\noindent We shall refer to a system ${\cal N}$ satisfying these constraints as an \emph{\nmsystem{}}.
Szemer\'edi and Petruska \cite{SzemerediPetruska} conjectured the following:

\begin{conjecture}\label{SPConjecture} Any \nmsystem{} satisfies $n \leq {m+2 \choose 2}$.
\end{conjecture}
Szemer\'edi and Petruska give a construction to show that this conjecture, if true, would be sharp.
Indeed it has been conjectured (by us and others) that this construction is the unique extremal structure for $m\geq 4$.
Gy\'arf\'as, Lehel, and Tuza \cite{GyarfasLehelTuza} proved $n\leq 2m^2+m$.  Later
Tuza \cite{Tuza2} proved the best known bound\footnote{The bound can be derived by 
by applying Theorem 5 from \cite{Tuza2} to the bound in the proof of Theorem 2 from \cite{Tuza2} and finishing with the bound from Proposition 7 in \cite{Tuza}.}, 
$n\leq \frac{3}{4}m^2+m+1$.  
These bounds use the machinery of $\tau$-critical hypergraphs.
Here we present an alternative approach that adapts an iterative decomposition process introduced by Szemer\'edi and Petruska and applies the skew version of Bollob\'as's theorem \cite{Bollobas} on the size of intersecting set pair systems. 
The skew version of Bollob\'as's theorem was first proven by Frankl \cite{Frankl} and also Kalai \cite{Kalai} 
(see also Theorem 5.6 of Babai and Frankl's book \cite{BabaiFrankl}).
While the bound we obtain, $n\leq m^2 + 6m + 2$ (Theorem \ref{maintheorem}), is weaker than Tuza's bound,
it is a proof-of-concept for an different method and a call to apply dimension bounds from linear algebra.

As noted by Gy\'arf\'as, Lehel, and Tuza \cite{GyarfasLehelTuza}, the Szemer\'edi and Petruska problem is equivalent
to determining the maximum order of a $\tau$-critical $3$-uniform hypergraph with transversal number $m$.  
They also determined that $O(m^{r-1})$ is the correct order of magnitude for the
maximum order of a $\tau$-critical $r$-uniform hypergraph with transversal number $m$.
Tuza \cite{Tuza2} proved the best known bounds.
The
methods presented offer an alternative approach and may yield improved bounds.  This is the focus of future research
which also will explain a connection to a conjecture of Lehel and Tuza (Problem $18$ of \cite{Tuza}) and a theorem of Hajnal \cite{Hajnal}.

Section \ref{decomposition} introduces notation and
recalls the process, introduced by Szemer\'edi and Petruska, to decompose \nmsystem{}s.
Section \ref{PairSelection} introduces a recursive procedure, based on this decomposition process,
to select special private pairs.  Section \ref{freepairs} defines a large subset of free private pairs chosen 
from this selection of special private pairs.
A skew $(2,m)$-system ultimately arises from this subset
of free private pairs in Section \ref{Skew}, where Theorem \ref{maintheorem}
is finally presented.  The last section includes remarks on possible improvements.


\section{The Decomposition Process} \label{decomposition}

We begin by giving definitions and recalling the process, introduced by Szemer\'edi and Petruska\footnote{We have endeavored to use the same notation introduced by Szemer\'edi and Petruska.  Important exceptions include that they use $n$ and $k$ to mean the quantities we refer to as $k$ and $\ell$, respectively.}, to decompose \nmsystem{}s.  
Much of the presentation in this section is lifted verbatim from their paper \cite{SzemerediPetruska}.
We assume $\ell\geq 4$ (Szemer\'edi and Petruska resolve the $\ell=2,3$ cases).
Let ${\cal N} = \{N_1,\ldots,N_\ell\}$, be an \nmsystem{}.
Define a collection of objects iteratively in {\em stages}, which are also called {\em times}, starting with stage $0$.
Set $\ell_0 = \ell$, ${\cal N}^{(0)} = {\cal N}$ and $N_i^{(0)} = N_i$.
For $i=1,\ldots,\ell_0$, fix a choice of vertex $x_i^{(0)} \in \bigcap_{j\neq i} N_j$.  By definition, $x_i^{(0)} \neq x_j^{(0)}$, for $i \neq j$.  
The set $A^{(0)} = \left\{ x_1^{(0)},\ldots,x_{\ell_0}^{(0)} \right\}$ is called the {\em kernel} at stage $0$; $x_1^{(0)},\ldots,x_{\ell_0}^{(0)}$ are the {\em kernel vertices} at stage $0$.

Assume at stage $j$ ($j \geq 0$) that $\ell_j\geq 4$, ${\cal N}^{(j)} = \left\{N_1^{(j)},\ldots,N_{\ell_j}^{(j)}\right\}$, and $A^{(j)} = \left\{ x_1^{(j)},\ldots,x_{\ell_j}^{(j)} \right\}$ are defined. Also assume that the minimal substructures of the ``remainder'' structure
$$
R^{(j)} = \left\{N_1^{(j)} - \bigcup_{i=0}^j A^{(i)},\ldots,N_{\ell_j}^{(j)}- \bigcup_{i=0}^j A^{(i)}\right\}
$$
satisfy (i).  We now explain the definition of $\ell_{j+1}$, ${\cal N}^{(j+1)}$, and $A^{(j+1)}$.  Consider substructures in $R^{(j)}$ that are minimal structures with respect to property (i).  Stop if there are no such substructures with more than three sets.  Otherwise, let
$$
{\cal N}^{(j+1)} = \left\{N_1^{(j+1)},\ldots,N_{\ell_{j+1}}^{(j+1)}\right\} \subset {\cal N}^{(j)}
$$
be chosen so that $\ell_{j+1} \geq 4$ and the corresponding remainders
$$\left\{N_1^{(j+1)} - \bigcup_{i=0}^{j} A^{(i)},\ldots,N_{\ell_{j+1}}^{(j+1)}- \bigcup_{i=0}^{j} A^{(i)}\right\}$$
form a substructure in $R^{(j)}$ that satisfies (i).
For $i=1,\ldots,\ell_{j+1}$, fix a choice of vertex 
$$x_i^{(j+1)} \in \bigcap_{r=1 \atop r\neq i}^{\ell_{j+1}} \left( N_r^{(j+1)} - \bigcup_{i=0}^{j} A^{(i)}\right).$$

By definition, $x_r^{(j+1)} \neq x_s^{(j+1)}$, for $1\leq r < s \leq \ell_{j+1}$.
The {\em kernel at stage $j+1$} is 
$$A^{(j+1)} = \left\{ x_1^{(j+1)},\ldots,x_{\ell_{j+1}}^{(j+1)} \right\}.$$
This process defines $\ell_{j}$, ${\cal N}^{(j)}$, and $A^{(j)}$ for $0 \leq j \leq t$, for some $t$.
Note that $t$ has been defined here as the length of the iterative process.

Because ${\cal N} = \left\{N_1,\ldots,N_{\ell_0}\right\}$ is an arbitrary enumeration of ${\cal N}$, we may assume that
$$
{\cal N}^{(j)} = \left\{N_1,\ldots,N_{\ell_{j}}\right\}, \mbox{ for } j = 0,\ldots,t.
$$
Define, for $i=1,\ldots,\ell_0$, the last {\em time} (or stage) that the truncation of $N_i$ appears 
in a substructure of this decomposition process, denoted $t_i$, as
$$
t_i = \max\left\{j : N_i \in {\cal N}^{(j)}\right\}.
$$
By definition, $t = t_1 \geq \cdots \geq t_{\ell_0} \geq 0$.

The next lemma gathers several properties of this iterative process.
\begin{lemma} \label{Observations} Some observations:
\begin{itemize}
\item[(a)] $A^{(j)}$ are pairwise disjoint, for $j=0,\ldots,t$.
\item[(b)] $\ell_j =\left|A^{(j)}\right| \geq 4$, for $j=0,\ldots,t$.
\item[(c)] $\left|N_i \cap A^{(r)} \right| = \ell_r - 1$, for $1 \leq i \leq \ell$, $0 \leq r \leq t_i$.
\item[(d)] $t < m$.
\item[(e)] $\ell_0 + \cdots + \ell_t \geq k - 2m = n - 3m$.
\end{itemize}
\end{lemma}
\begin{proof} Properties (a) -- (c) are immediate from the definition.  Properties (d) and (e) are respectively Lemma 5 and Lemma 6 of \cite{SzemerediPetruska}.
\end{proof}

Define $A = \bigcup_{i=0}^{t} A^{(i)}$; it is the set of kernel vertices.  Let $G = V - A$ denote the {\em garbage} vertices; that is, the vertices remaining after the aforementioned kernel-defining decomposition process terminates.

\section{Selection of private pairs} \label{PairSelection}

In this section we define a process to select private pairs.  Much of the beginning of this
is a review of results from the paper by Szemer\'edi and Petruska \cite{SzemerediPetruska}.

A pair of elements $p \subset N_i$ is {\em single-covered with respect to ${\cal N}^{(j)}$}, for some $j$ satisfying $0 \leq j \leq t_i$, if $N_i$ is the only set in ${\cal N}^{(j)}$ that contains $p$ as a subset.
If $p \subset N_i$ is single-covered with respect to ${\cal N}^{(j)}$,
then it is a {\em private pair for $N_i$ at time $j$}, or simply a {\em private pair}.  Observe that if a pair
is private for $N_i$ at time $j$, then it remains a private pair for $N_i$ until (and including) time $t_i$.
A pair that is contained in at least two sets in ${\cal N}^{(j)}$ is called a
{\em double-covered} pair (at time $j$).

The following lemma is a reformulation of Lemma 7 from 
\cite{SzemerediPetruska}.
\begin{lemma} \label{YLemma} For all $j=0,\ldots,t$,
\begin{itemize}
\item[(a)] Every pair from $A$ is double-covered by ${\cal N}^{(j)}$.  
\item[(b)] For all $N_i \in {\cal N}^{(j)}$ and any subset $Y \subseteq N_i$ such that $|Y| = j$,
there exists a pair in $N_i - Y$ 
that is single-covered with respect to ${\cal N}^{(j)}$.
\end{itemize}
\end{lemma}
\begin{proof} Properties (a) and (b) are, respectively, the proof of Lemma 7 part (b) and the proof of Lemma 7 part (a) of \cite{SzemerediPetruska}.
\end{proof}

An important consequence of Lemma \ref{YLemma} (a) is the following.
\begin{corollary} \label{GarbageVertexInPrivatePair} Any private pair must contain at least one vertex from $G$.
\end{corollary}

Next we describe a process to select a collection of private pairs for each $N_i$.
For each $i \in \{1,\ldots,\ell\}$, we define, by induction on time, a 
set $P_i = \left\{p_i^{(j)} : 0 \leq j \leq t_i \right\}$ of $t_i+1$ private pairs for $N_i$.  
The pair $p_i^{(j)}$ will be chosen from among the private pairs for $N_i$ at time $j$.  
By Corollary \ref{GarbageVertexInPrivatePair} the pair $p_i^{(j)}$ contains at least one vertex from $G$.
Choose $g_i^{(j)} \in p_i^{(j)} \cap G$; it is the {\em anchor} of $p_i^{(j)}$.
The other element of $p_i^{(j)}$ is the {\em non-anchor} of $p_i^{(j)}$; it is denoted $u_i^{(j)}$.
Naturally it is possible that the non-anchor is also an element of $G$, but we shall distinguish $g_i^{(j)}$ as the anchor.
Auxiliary sets $P_i^{(j)}$ and $G_i^{(j)}$ will also be defined; 
they are the initial segments of the private pairs and anchors for the private pairs selected for $N_i$.  

Initially, for $i \in \{1,\ldots,\ell\}$, 
let $p_i^{(0)}= \left\{u_i^{(0)}, g_i^{(0)}\right\} $ be a private pair for $N_i$ at time zero.
Such a private pair exists by an application of Lemma \ref{YLemma} part (b) in which $Y=\varnothing$.
Set $G_i^{(0)} = \left\{ g_i^{(0)} \right\}$ and $P_i^{(0)} = \left\{ p_i^{(0)} \right\}$.

For $j>0$ and $i \in \{1,\ldots,\ell\}$, assume that the sets $P_i^{(j-1)}$ and $G_i^{(j-1)}$ have already been defined.  Also assume that a private pair $p_i^{(j)}=\left\{u_i^{(j)}, g_i^{(j)}\right\}$ has already been chosen for each $N_i$ with $j \leq t_i$.
Now define 
\stack{$G_i^{(j)}$}{$=$}{$G_i^{(j-1)} \cup \left\{g_i^{(j)}\right\}$}{if $j \leq t_i$}{$G_i^{(j-1)}$}{if $j > t_i$,}
and similarly define,
\stack{$P_i^{(j)}$}{$=$}{$P_i^{(j-1)} \cup \left\{p_i^{(j)}\right\}$}{if $j \leq t_i$}{$P_i^{(j-1)}$}{if $j > t_i$.}
This definition yields $P_i^{(j)}=\left\{p_i^{(0)},\ldots,p_i^{(j)}\right\}$ and $G_i^{(j)}=\left\{g_i^{(0)},\ldots,g_i^{(j)}\right\}$, for all $0 \leq j \leq t_i$.  In particular,
note that $\left|P_i^{(j)}\right| = \left|G_i^{(j)}\right| = j+1$, for all $0 \leq j \leq t_i$.
Also
$\bigcup_{h=0}^\ell P_h^{(j)}$ represents the set of private pairs defined up through time $j$.

To complete the iterative process, it 
remains to describe how to select a private pair $p_i^{(j)}$, for all $i \in \{1,\ldots,\ell_j\}$.  
Apply Lemma \ref{YLemma} part (b) with $Y=G_i^{(j-1)}$ to produce a pair $p_i^{(j)}= \left\{u_i^{(j)}, g_i^{(j)}\right\}$ private for $N_i^{(j)}$ satisfying 
$p_i^{(j)} \cap Y = \varnothing$.  In particular,
$g_i^{(j)} \in G$ and $g_i^{(j)} \notin Y = G_i^{(j-1)}=\left\{g_i^{(0)},\ldots,g_i^{(j-1)}\right\}.$

If the non-anchor of $p_i^{(j)}$ is in $\bigcup_{s=0}^{j-1} A^{(s)}$, say $u_i^{(j)} = x_a^{(b)}$, for some $b < j$, then replace $u_i^{(j)}$ with either $x_a^{(j)}$, if $j\leq t_a$, or $x_a^{(t_a)}$ otherwise.  Observe that after this replacement the new pair is still private to $N_i$ at time $j$.  In other words, if 
$u_i^{(j)} \in \bigcup_{s=0}^{j-1} A^{(s)}$, then 
$u_i^{(j)} = x_a^{(t_a)}$, for some $a \in \{1,\ldots,\ell\}$.

Finally, suppose that $u_i^{(j)} \in  A^{(j)}$, say $u_i^{(j)} = x_a^{(j)}$, for some $a \in \{1,\ldots,\ell\}$.  If $j < t_a$ and $j < t_i$ then
set $u_i^{(j)} = x_a^{(j+1)}$.  After this replacement the new pair is still private to $N_i$ at time $j$.  In other words, if 
$u_i^{(j)} = x_a^{(j)}$, for some $a \in \{1,\ldots,\ell\}$, then $j \in \{t_a,t_i\}$.

Let $P^{(j)} = \left\{p_i^{(j)} : 1 \leq i \leq \ell_j \right\}$ denote the private pairs defined by this process at time $j$ and let $P_i = \left\{p_i^{(j)} : 0 \leq j \leq t_i \right\}$ be the set of $t_i+1$ private pairs defined for $N_i$ by this process.  The collection of all 
selected pairs is defined as
$$
P = \bigcup_{j=0}^t P^{(j)} =\bigcup_{i=0}^\ell P_i.
$$ 
\noindent
The following lemma gathers observations about private pairs in $P$.

\begin{lemma} \label{importantobservations} The pairs in $P$ satisfy: 
\begin{itemize}
\item[(a)] $P^{(j_1)} \cap P^{(j_2)} = \varnothing$, for $0 \leq j_1 < j_2 \leq t$.
\item[(b)] Any pair in $\bigcup_{s=0}^j P^{(s)}$ is at most single-covered by ${\cal N}^{(j)}$, for $j = 0,\ldots,t$.
\item[(c)] $\left|P^{(j)}\right| = \ell_j$ and every $\left|P_i \cap P^{(j)}\right| = 1$, for all $0 \leq j \leq t$ and $1 \leq i \leq \ell_j$.
\item[(d)] $|P| = \left| \bigcup_{j=0}^t P^{(j)} \right| = \sum_{j=0}^t \ell_j \geq n - 3m$.
\item[(e)] If 
$u_i^{(j)} \in \bigcup_{s=0}^{j-1} A^{(s)}$, then 
$u_i^{(j)} = x_a^{(t_a)}$, for some $a \in \{1,\ldots,\ell\}$.
\item[(f)] If 
$u_i^{(j)} = x_a^{(j)}$, for some $a \in \{1,\ldots,\ell\}$, then $j \in \{t_a,t_i\}$.
\end{itemize}
\end{lemma}
\begin{proof} 
Parts (a), (b), and (c) , respectively, follow from the same arguments given to prove parts (*), (**), and (***) of Lemma 8 of \cite{SzemerediPetruska}.
Part (d) is essentially a consequence of the arguments given to prove Lemma \ref{Observations}(e).
Part (e) and (f) reiterate the observations in the paragraphs defining the selection of the private pairs in $P$.
\end{proof}

\section{Free pairs} \label{freepairs}

In this section we define a special subset of private pairs in $P$ that is used in Section \ref{Skew} to define a large skew $(2,m)$-system.
Recall that $A$ is the set of kernel vertices.  Every element of $A$ has the form
$x_i^{(j)}$, where $1 \leq i \leq \ell$ and $0 \leq j \leq t_i$.  

Create a digraph $D$ on the vertex set $A$ in which there is an
arc from $x_r^{(s)}$ to $x_i^{(j)}$ if $u_r^{(s)} = x_i^{(j)}$ and $j \neq t_i$.

\begin{lemma} \label{Arcs} If there is an arc in $D$ from $x_r^{(s)}$ to $x_i^{(j)}$ then
\begin{itemize}
\item[(a)] $s \leq j$, and
\item[(b)] if $s = j$, then $s = t_r$.
\end{itemize}
\end{lemma}
\begin{proof} Because $j \neq t_i$, Lemma \ref{importantobservations}(e) implies that
$s \leq j$.  If $s = j$, then Lemma \ref{importantobservations}(f) guarantees
that $j \in \{t_r,t_i\}$.  Consequently $t_r = j = s$.
\end{proof}

\begin{lemma}  The digraph $D$ is acyclic and has out-degree at most one.
\end{lemma}
\begin{proof} Lemma \ref{importantobservations} part (c) guarantees that the out-degree of vertex  $x_r^{(s)}$ in $D$ is at most one.  Suppose, to the contrary, that there are arcs forming a directed cycle:
$$
x_{r_1}^{(s_1)} \rightarrow x_{r_2}^{(s_2)} \rightarrow \cdots \rightarrow x_{r_h}^{(s_h)} \rightarrow x_{r_1}^{(s_1)}. 
$$
Lemma \ref{Arcs}(a) yields $s_1 \leq s_2 \leq \cdots \leq s_h \leq s_1$,
so $s_1 = s_i$, for all $1 \leq i \leq h$.  Lemma \ref{Arcs}(b) then implies that 
$t_{r_i} = s_1$, for all $1 \leq i \leq h$.  
But the arc from $x_{r_h}^{(s_h)}$ to $x_{r_1}^{(s_1)}$ requires $s_1 \neq t_{r_1}$.
\end{proof}

The graph obtained from $D$ by removing direction on the arcs is a forest that contains a maximum independent set of vertices; call it $F$.
Because forests are $2$-colorable, it follows that $|F| \geq |A|/2$.  Define free pairs in $P$ this way:
a pair $p_r^{(s)} \in P$ is {\em free} if $x_r^{(s)} \in F$.

The most important consequence of Lemma \ref{importantobservations} is a lower bound on the number of free pairs.
\begin{corollary}  \label{NumberFreePairs} The number of free pairs in $P$ is at least $\frac{n - 3m}{2}$.
\end{corollary}
\begin{proof} By part (d) of Lemma \ref{importantobservations}, $|A| \geq n - 3m$.
Because$|F| \geq |A|/2$, the result follows.
\end{proof}

\section{A skew system} \label{Skew}

In this section we apply the following theorem, first proven by Frankl \cite{Frankl}; it is the skew version of a theorem due to Bollob\'as \cite{Bollobas}.  This theorem is also 
presented in the book by Babai and Frankl (\cite{BabaiFrankl}, pages $94$--$95$).

\begin{theorem} \label{SkewBollobas} (Bollob\'as's Theorem - Skew Version)
If $A_1,\ldots,A_h$ are $r$-element sets and $B_1,\ldots,B_h$ are $s$-element sets such that
\begin{itemize}
\item[(a)] $A_i$ and $B_i$ are disjoint for $i=1,\ldots,h$,
\item[(b)] $A_i$ and $B_j$ intersect whenever $1 \leq i < j \leq h$
\end{itemize}
then $h \leq {r+s \choose r}$.
\end{theorem}

A system of sets, $\{(A_i,B_i)\}_{i=1}^h$, satisfying the hypotheses of Theorem \ref{SkewBollobas} is called a {\em skew intersecting set pair $(r,s)$-system}; abbreviate this to {\em skew $(r,s)$-system}.

The goal in this section is to apply  Theorem \ref{SkewBollobas} to a skew $(2,m)$-system derived from the free pairs in $P$.
First use all of the pairs in $P$ to define, iteratively, a collection of $m$-sets this way.
To each $N_i$ associate $t_i+1$ $m$-sets denoted $M_i^{(0)},\ldots,M_i^{(t_i)}$.
At stage $0$, set $M_i^{(0)} = \overline{N_i}$, for all $i=1,\ldots,\ell_0$.
For $i=1,\ldots,\ell$ and $j = 1,\ldots,t_i$, recursively define
\stack{$M_i^{(j)}$}{$=$}{$M_i^{(j-1)} - \left\{x_i^{(j-1)}\right\} + \left\{g_i^{(j-1)}\right\}$}{if $p_i^{(j-1)}$ is free}{$M_i^{(j-1)}$}{if $p_i^{(j-1)}$ is not free}
Note that, because $\left|M_i^{(0)}\right|=m$, it follows that $\left|M_i^{(j)}\right| = m$, for all 
$i=1,\ldots,\ell$ and $j = 1,\ldots,t_i$.  
Also observe that this
recursive process will never remove $x_i^{(t_i)}$ from $M_i^{(0)}$
because the process halts at stage $t_i$.

Now define the system
$${\cal F} = \left\{ \left(p_i^{(j)},M_i^{(j)}\right): p_i^{(j)} \in P \mbox{ is free}\right\},$$
where ${\cal F}$ is ordered linearly and chronologically via lexicographical order:
$$\Bigl(p_r^{(s)},M_r^{(s)}\Bigr) < \left(p_i^{(j)},M_i^{(j)}\right) \qquad \Longleftrightarrow \qquad \mbox{($s < j$) or ($s=j$ and $r < i$}).$$

\begin{theorem} \label{SkewSystem} ${\cal F}$ is a skew $(2,m)$-system.
\end{theorem}
\begin{proof} Clearly $\left|p_i^{(j)}\right| = 2$ and $\left|M_i^{(j)}\right| = m$, for
all $\left(p_i^{(j)},M_i^{(j)}\right) \in {\cal F}$.  Because $p_i^{(j)}$ is private to $N_i$ at time $j$,
it follows that $p_i^{(j)} \subset N_i$; so, $p_i^{(j)} \cap M_i^{(0)} = \varnothing$.  Observe 
that the anchor, $g_i^{(j)}$, for $p_i^{(j)}$ can not be added to $M_i^{(j)}$ by the recursive process generating the $m$-sets because $g_i^{(j)}$ is added at time $j+1$.  
Moreover, by the iterative choice of private pairs, the pair
$p_i^{(j)}$ was chosen to be disjoint from $\left\{g_i^{(0)},\ldots, g_i^{(j-1)}\right\}$, implying that the
non-anchor of $p_i^{(j)}$ is never added to $M_i^{(0)}$ in the process to produce $M_i^{(j)}$.
Therefore, $p_i^{(j)} \cap M_i^{(j)} = \varnothing$, showing that hypothesis (a) is satisfied in Theorem \ref{SkewBollobas}.

Now we prove the system satisfies hypothesis (b). 
Suppose $\left(p_r^{(s)},M_r^{(s)}\right), \left(p_i^{(j)},M_i^{(j)}\right) \in {\cal F}$ and $\left(p_r^{(s)},M_r^{(s)}\right) < \left(p_i^{(j)},M_i^{(j)}\right)$.  We must prove $p_r^{(s)} \cap M_i^{(j)} \neq \varnothing$.
If $r = i$, then $s < j$ so $g_r^{(s)}$ is in $M_i^{(s+1)}$ (and therefore $M_i^{(j)}$) because $p_r^{(s)}$ is free.
Consequently we may assume $r \neq i$.

If $g_r^{(s)} \in M_i^{(j)}$, then $p_r^{(s)} \cap M_i^{(j)} \neq \varnothing$.   So we may assume
$g_r^{(s)} \notin M_i^{(j)}$.  
Elements from $G$ are
only added to $M_i^{(0)}$ to get to $M_i^{(j)}$, so $g_r^{(s)} \notin M_i^{(0)} = \overline{N_i}$ which implies $g_r^{(s)} \in N_i$.  
Since $s\leq j$ and $p_r^{(s)}$ is private to $N_r$ at time $s$, we conclude that 
$u_r^{(s)} \notin N_i$.  So $u_r^{(s)} \in \overline{N_i} = M_i^{(0)}$. 
 
If $u_r^{(s)} \neq x_i^{(\alpha)}$, for some $0 \leq \alpha \leq j-1$,
then $u_r^{(s)} \in M_i^{(j)}$ meaning $p_r^{(s)} \cap M_i^{(j)} \neq \varnothing$.

So we may assume that 
$u_r^{(s)} = x_i^{(\alpha)}$, for some $0 \leq \alpha \leq j-1$.
Because $\alpha < j \leq t_i$, note that $\alpha \neq t_i$.
By definition there is arc in the digraph $D$ from $x_r^{(s)}$ to $x_i^{(\alpha)}$.
But $p_r^{(s)}$ is free, so $x_r^{(s)} \in F$.  Since $F$ is an independent set in $D$,
it follows that $x_i^{(\alpha)} \not\in F$.
Accordingly $x_i^{(\alpha)}$ is never removed from 
$M_i^{(0)}$ in the production of $M_i^{(j)}$.  Therefore $x_i^{(\alpha)} = u_r^{(s)} \in p_r^{(s)} \cap M_i^{(j)}$.
\end{proof}


Now we are ready to state the main theorem of the paper.

\begin{theorem} \label{maintheorem}
Any \nmsystem{} satisfies $n\leq m^2 + 6m + 2$.
\end{theorem}
\begin{proof} Corollary \ref{NumberFreePairs} proves $\frac{n - 3m}{2} \leq |{\cal F}|$.
Combining Theorem \ref{SkewBollobas} with Theorem \ref{SkewSystem} yields
$
|{\cal F}| \leq {m+2 \choose 2};
$
therefore $n\leq m^2 + 6m + 2$.
\end{proof}

It is possible to reduce the upper bound in Theorem \ref{maintheorem}.
Natural reductions can be achieved in two ways: increasing the bound on $|P|$ that appears in 
Lemma \ref{importantobservations}(d) and enlarging the set of pairs used to define ${\cal F}$.
However the improvements that we have found do not reduce the leading term of the bound and so for simplicity's sake we have opted
to omit them.  

More tantalizingly hopeful is a linear algebra approach (using a dimension argument similar to Lov\'{a}sz's proof 
of Bollob\'as's theorem as presented in the book by Babai and Frankl \cite{BabaiFrankl}) to prove
Conjecture \ref{SPConjecture}.  Small computations confirm the linear independence of appropriately chosen
homogeneous polynomials of degree two in $m+1$ variables associated with carefully selected private pairs for each
$N_i$ and each
vertex from $G$.  
Unfortunately a general proof of the linear independence of these polynomials has eluded us.
The argument presented here essentially uses the skew version of Bollob\'as's theorem to verify the linear
independence of a large number of these polynomials.  Numerous small extremal structures and the unwieldy
form of the general conjecture (Problem $18$ of \cite{Tuza}) also pose serious obstacles for larger uniformity ($r > 3$).



\section*{Acknowledgments}
We thank Jen\H{o} Lehel for helpful remarks and discussions.
We are also very grateful to Zsolt Tuza for bringing our attention to reference \cite{Tuza2} and alerting
us to the upper bound $\frac{3}{4}m^2+m+1$ for the 
maximum order of a $\tau$-critical $3$-uniform hypergraph with transversal number $m$.


\end{document}